\documentclass[12pt]{amsart}

\usepackage{mathrsfs}
\usepackage{amsfonts}
\usepackage{amssymb}
\usepackage{amsfonts, amscd, amsmath, mathrsfs, amssymb, amsthm, amsxtra, bbding, epsfig, graphicx, latexsym, url, mathbbol, bbold}
\usepackage[papersize={8.2in,11in},textwidth=16.2cm,textheight=22cm,centering]{geometry}
\usepackage{enumerate}
\usepackage{caption} 

\usepackage{graphicx}
\usepackage{tikz}
\usepackage{float}

\usepackage{xcolor}
\definecolor{cite}{rgb}{0.00,0.60,1.00}
\definecolor{url}{rgb}{1.00,0.10,0.80}
\definecolor{link}{rgb}{0.00,0.00,1.00}
\usepackage[colorlinks,linkcolor=link,urlcolor=url,citecolor=cite,pagebackref,breaklinks]{hyperref}

\hypersetup{
	pdfstartpage=1,
	pdfstartview=FitH}

\DeclareFontFamily{U}{mathx}{\hyphenchar\font45}
\DeclareFontShape{U}{mathx}{m}{n}{
	<5> <6> <7> <8> <9> <10>
	<10.95> <12> <14.4> <17.28> <20.74> <24.88>
	mathx10
}{}
\DeclareSymbolFont{mathx}{U}{mathx}{m}{n}
\DeclareMathAccent{\widecheck}{\mathalpha}{mathx}{"71}

\numberwithin{equation}{section}

\allowdisplaybreaks

\newtheorem{theorem}{Theorem}[section]
\newtheorem{lemma}{Lemma}[section]
\newtheorem{problem}{Problem}[section]

\newtheorem{proposition}[theorem]{Proposition}
\newtheorem{example}{Example}[section]
\newtheorem{corollary}{Corollary}[section]

\makeatletter
\newcounter{roem}
\renewcommand{\theroem}{\Roman{roem}}

\newcommand{\c@org@eq}{}
\let\c@org@eq\c@equation
\newcommand{\org@theeq}{}
\let\org@theeq\theequation

\newcommand{\setroem}{
	\let\c@equation\c@roem
	\let\theequation\theroem}

\newcommand{\setarab}{
	\let\c@equation\c@org@eq
	\let\theequation\org@theeq}
\makeatother

\newtheorem{remark}{\bf Remark}

\newcommand{\N}{\mathbb{N}}

\newcommand{\E}{\mathbb{E}}

\newcommand{\bP}{\mathbf{P}}

\def\le{\leqslant}
\def\leq{\leqslant}
\def\ge{\geqslant}
\def\geq{\geqslant}

\usepackage{graphicx}
\usepackage{tikz}

\begin{document}

	\title[Multiple representations in difference sets]{Multiple representations in difference sets under polynomial density constraints}

	\author
[Y. Ding, S. Kiss, W. Niu, P. P. Pach and C. S\'andor] {Yuchen Ding, S\'andor Z. Kiss, Wei Niu, P\'eter P\'al Pach and Csaba S\'andor}
	
\address{(Yuchen Ding$^{1,2}$) $^1$School of Mathematics,  Yangzhou University, Yangzhou 225002, People's Republic of China}
\address{$^2$HUN-REN Alfr\'ed R\'enyi Institute of Mathematics, Budapest, Pf. 127, H-1364 Hungary}
\email{ycding@yzu.edu.cn}

\address{(S\'andor Kiss)  Department of Algebra and Geometry, Institute of Mathemetics, Budapest University of Technology and Economics, M\H{u}egyetem rkp. 3., H-1111 Budapest, Hungary; 
\newline \hspace*{4mm}
HUN-REN Alfr\'ed R\'enyi Institute of Mathematics, Re\'altanoda utca 13--15., H-1053 Budapest, Hungary}
\email{kiss.sandor@ttk.bme.hu}

\address{(Wei Niu$^{1,2}$) $^1$School of Mathematics and Statistics, Xi’an Jiaotong University, 28 Xianning
West Road, Beilin District, Xi’an, Shaanxi, China}
\address{$^2$Department of Computer Science and Information Theory,
Budapest University of Technology and Economics,
M\H{u}egyetem rakpart 3,
H-1111 Budapest,
Hungary}
\email{wei.niu@stu.xjtu.edu.cn}

\address{(P\'eter P\'al Pach)  HUN-REN Alfr\'ed R\'enyi Institute of Mathematics, Re\'altanoda utca 13--15., H-1053 Budapest,  Hungary; \newline \hspace*{4mm}
MTA--HUN-REN RI Lend\"ulet ``Momentum'' Arithmetic Combinatorics Research Group, Re\'altanoda utca 13--15., H-1053 Budapest,  Hungary; \newline \hspace*{4mm}
Department of Computer Science and Information Theory, Budapest University of Technology and Economics, M\H{u}egyetem rkp. 3., H-1111 Budapest, Hungary; \newline \hspace*{4mm}
Extremal Combinatorics and Probability Group (ECOPRO), Institute for Basic Science (IBS), Daejeon, South Korea.}
\email{pachpp@renyi.hu}

\address{(Csaba S\'andor)  Department of Stochastics, Institute of Mathemetics, Budapest University of Technology and Economics, M\H{u}egyetem rkp. 3., H-1111 Budapest, Hungary; \newline \hspace*{4mm}
HUN-REN Alfr\'ed R\'enyi Institute of Mathematics, Re\'altanoda utca 13--15., H-1053 Budapest,  Hungary; \newline \hspace*{4mm} 
MTA--HUN-REN RI Lend\"ulet ``Momentum'' Arithmetic Combinatorics Research Group, Re\'altanoda utca 13--15., H-1053 Budapest,  Hungary}
\email{sandor.csaba@ttk.bme.hu}

\begin{abstract}
Motivated by the classical study of additive representation functions and the combinatorial phenomenon of sum-difference asymmetry, we investigate the subtractive representation function $d_{\mathcal{B},\mathcal{A}}(n)$, which counts the number of solutions of the equation $n = b - a$ with $a\in \mathcal{A}$ and $b\in \mathcal{B}$. Specifically, we study lower bounds and the asymptotic behavior of $d_{\mathcal{B},\mathcal{A}}(n)$ under polynomial density constraints of orders $x^\alpha$ and $x^{1-\alpha}$, where $0 < \alpha < 1$. We first establish a general lower bound for the first moment of $d_{\mathcal{B},\mathcal{A}}(n)$. Furthermore, by employing a recursive construction, we prove that this lower bound is essentially optimal under certain regularity assumptions on the counting functions. Finally, utilizing probabilistic methods, we show the existence of sets for which $d_{\mathcal{B},\mathcal{A}}(n)$ has a logarithmic point-wise upper bound. 

{\it 2010 Mathematics Subject Classification:} 11B34

{\it Keywords and phrases:} additive number theory; difference set; probabilistic method  

\end{abstract}

	\maketitle

	\setcounter{tocdepth}{1}
	
	
	\section{Introduction}

Let $\mathbb{N}$ be the set of natural numbers and $\mathcal{A}, \mathcal{B}$ be two subsets of $\mathbb{N}$. The set $\mathcal{B}$ is said to be an additive complement of $\mathcal{A}$ if all sufficiently large integers can be written as the sum of an element of $\mathcal{A}$ and an element of $\mathcal{B}$. The study of additive complements and their associated representation functions, namely $r_{\mathcal{A},\mathcal{B}}(n) = \#\{(a,b) \in \mathcal{A} \times \mathcal{B} : n = a+b\}$, is a classic topic in additive number theory. Clearly, $$\sum_{n\le x}r_{\mathcal{A},\mathcal{B}}(n)\ge x-C$$ for some constant $C\geq 0$,  if the set $B\subseteq \mathbb{N}$ is the additive complement of the set $A\subseteq \mathbb{N}$. 

For any $\mathcal{S} \subseteq \mathbb{N}$, let $\mathcal{S}(x)$ denote the number of elements of the set $\mathcal{S}$ up to $x$. Obviously, 
if $c_{1}x^\alpha < \mathcal{A}(x) < c_{2}x^{\alpha}$ and $c_{3}x^{1-\alpha} < \mathcal{B}(x) < c_{4}x^{1-\alpha}$, where $c_{1}, c_{2}, c_{3}, c_{4}$ are positive constants and $0<\alpha<1$, then
\[
c_{5}x \le \mathcal{A}\left(\frac{x}{2}\right)\mathcal{B}\left(\frac{x}{2}\right) \le \sum_{n\le x}r_{\mathcal{A},\mathcal{B}}(n) \le \mathcal{A}\left(x\right)\mathcal{B}\left(x\right) \le c_{6}x, 
\]
where $c_{5}, c_{6}$ are positive. For a fixed $\ell \ge 2$, consider the set of perfect powers $\mathcal{S}_{\ell} = \{1^{\ell}, 2^{\ell}, \dots\}$, and let $\mathcal{W}_{\ell}$ be the additive complement of $\mathcal{S}_{\ell}$.

Recently, motivated by Cilleruelo's conjecture \cite{Ci93} on the additive completion of $\ell$-th powers and related problems asked by Ruzsa \cite{Ru01} and Green (see \cite{CF17}), extensive research has been devoted to the exact growth rate of $r_{\mathcal{A},\mathcal{B}}(n)$ under specific algebraic structures or density constraints (see, for instance, \cite{DSZ26} and the references therein). 

In particular, Cilleruelo conjectured that 
\[
\sum_{n\le x}r_{\mathcal{S}_{\ell},\mathcal{W}_{\ell}}(n)-x \gg x.
\]
Recently for $\ell \ge 2$, it was proved in \cite{DSZ26} that
\[
\sum_{n\le x}r_{\mathcal{S}_{\ell},\mathcal{W}_{\ell}}(n) -x\gg x^{1-\frac{1}{\ell}}
\]
Furthermore, for $\ell = 2$,
\[
\sum_{n\le x}r_{\mathcal{S}_{2},\mathcal{W}_{2}}(n) -x\gg x^{\frac{3}{4}-\varepsilon}.
\]
The proofs are based on the following ideas. Let $K$ be a positive integer and let $N$ be a large integer. For $1 \le k \le K$, the authors define
\[
\mathcal{S}_{\ell,k,N} = \mathcal{S}_{\ell} \cap \left[\frac{k-1}{2K}N, \frac{k}{2K}N \right)
\]
and
\[
\mathcal{W}_{\ell,k,N} = \mathcal{W}_{\ell} \cap \left[\frac{k-1}{2K}N, \frac{k}{2K}N \right).
\]
In the proofs it was shown that 
\[
\sum_{1\le k\le K}|\mathcal{S}_{\ell,k,N}|\cdot |\mathcal{W}_{\ell,k,N}| \ge \frac{N\log K}{K}\left(1-\frac{1}{\ell}\right)\frac{1}{\ell}. 
\]
Introducing the notation
\[
d_{\mathcal{B},\mathcal{A}}(n) = \#\{(b,a)\in\mathcal{B}\times\mathcal{A}:n=b-a\},
\]
one has
\[
\sum_{-\frac{N}{2K}\le n \le \frac{N}{2K}}d_{\mathcal{S}_{\ell}, \mathcal{W}_{\ell}}(n) \gg \frac{N\log K}{K}.
\]
It follows that 
\[
\frac{\sum_{|n|\le x}d_{\mathcal{S}_{\ell}, \mathcal{W}_{\ell}}(n)}{x} \rightarrow \infty
\]
as $x\rightarrow \infty$. Based on these observations, in this paper we study the order of magnitude of 
\[
\sum_{0 < n \le x}d_{\mathcal{B},\mathcal{A}}(n),
\]
where $\mathcal{A}(x) \asymp x^{\alpha}$ or $\mathcal{A}(x) = (c+o(1))x^{\alpha}$ and 
$\mathcal{B}(x) \asymp x^{1-\alpha}$ as $x\rightarrow \infty$.


The function $d_{\mathcal{B},\mathcal{A}}(n)$ is already appeared in \cite{Nat07a} and \cite{NOORS} in a more general form; in these papers the representation functions associated to linear forms with integral coefficients were considered. 
In this article, we study lower bounds and the asymptotic behavior of $d_{\mathcal{B},\mathcal{A}}(n)$ under natural polynomial density constraints of orders $x^\alpha$ and $x^{1-\alpha}$, where $0 < \alpha < 1$.

From now on, let $0<\alpha<1$ be a given number. 
Let $f(x)$ be a function with nonnegative values. 
We use the standard Landau notation $O(\cdot)$, $\ll$, $\gg$, and $\asymp$; the implicit constants may depend on $\alpha$ and on the constants in the hypotheses, but not on $x$ or $n$. In particular, we write $\mathcal{S}(x)\asymp f(x)$ if there are two positive constants $c_1, c_2$ (not depending on $x$) such that
$c_1f(x)\le \mathcal{S}(x)\le c_2f(x)$. If $f(x)$ and $g(x)$ are real-valued functions, we write $f(x) = O(g(x))$ if there exists a constant $C > 0$ such that $|f(x)|\le C|g(x)|$ and $f(x)=o(g(x))$ with $g(x)\neq 0$ if
$
\lim_{x\to\infty}f(x)/g(x)=0
$. We also denote $f(x) = O(g(x))$ by $f(x) \ll g(x)$. The formula $f(x) \sim g(x)$ as $x \rightarrow \infty$ means 
\[
\lim_{x\rightarrow \infty}\frac{f(x)}{g(x)} = 1.
\]

We first give the following result.

	\begin{proposition}\label{thm:lb:ABasy}
		For any subsets $\mathcal{A}, \mathcal{B}\subset \N$ with $\mathcal{A}(x)\asymp x^{\alpha}$ and $\mathcal{B}(x)\asymp x^{1-\alpha}$, we have
		\[
		\sum_{n\leq x}d_{\mathcal{B},\mathcal{A}}(n) \gg x.
		\]
		\end{proposition}

The upper-bound hypothesis on $\mathcal{B}(x)$ in Proposition \ref{thm:lb:ABasy} is necessary; see Remark \ref{remark:upper_bound} in Section 2. The bound in Proposition \ref{thm:lb:ABasy} is sharp in general; see Example \ref{example:1}. That example does not have $\mathcal{A}(x)$ asymptotic to a constant multiple of $x^\alpha$, so we next consider what additional lower bounds follow when the counting function of $\mathcal{A}$ has a prescribed main term.

\begin{theorem}\label{thm:2}
Let \(f(x)\) be a monotone increasing function with $f(x)\le x^{\alpha }$ and \(f(x)\to\infty\) as \(x\to\infty\). Let also $\mathcal{A}, \mathcal{B}\subset \N$ be two subsets that satisfy
$$
\mathcal{A}(x)=cx^{\alpha}+O\left(\frac{x^\alpha}{f(x)}\right)
$$ 
and $\mathcal{B}(x)\asymp x^{1-\alpha}$, where $c>0$ is a given number.
Then we have
	\[
	\sum_{n\leq x}d_{\mathcal{B},\mathcal{A}}(n)\gg x\log f(x).
	\]
    Moreover, we also have
	$$
\limsup_{n\to\infty}\frac{d_{\mathcal{B},\mathcal{A}}(n)}{\log f(n)}>0.
$$
\end{theorem}

As several consequences of Theorem \ref{thm:2}, we mention the following corollaries. 

	\begin{corollary}\label{corollary-1}
    Let $c>0$ be a given number.
	Let $\mathcal{A}, \mathcal{B}\subset \N$ be two subsets that satisfy
    $\mathcal{A}(x)\sim cx^\alpha$ as $x\to\infty$ and  $\mathcal{B}(x)\asymp x^{1-\alpha}$. Then as $x\to\infty$ we have
	\[
\frac{1}{x}\sum_{n\leq x}d_{\mathcal{B},\mathcal{A}}(n)\to \infty.
	\]

\end{corollary}
Corollary \ref{corollary-1} shows that the lower bound of Proposition \ref{thm:lb:ABasy} can be significantly improved if $\mathcal{A}(x)$ has an asymptotic equivalent. 

\begin{corollary}\label{corollary-3}
    Let $c>0$ be a given number.
	Let $\mathcal{A}, \mathcal{B}\subset \N$ be two subsets that satisfy
    $\mathcal{A}(x)\sim cx^\alpha$ as $x\to\infty$ and  $\mathcal{B}(x)\asymp x^{1-\alpha}$. Then the sequence $d_{\mathcal{B},\mathcal{A}}(n)$ is unbounded.
\end{corollary}

\begin{corollary}\label{corollary-2}
Let $0<\beta<\alpha$ be two given numbers. If $\mathcal{A}, \mathcal{B}\subset \N$ are two subsets that satisfy
$\mathcal{A}(x)=cx^{\alpha}+O\big(x^{\beta}\big)$ and $\mathcal{B}(x)\asymp x^{1-\alpha}$, then we have
	\[
	\sum_{n\leq x}d_{\mathcal{B},\mathcal{A}}(n)\gg x\log x.
	\]
\end{corollary}

Our next theorem shows that Theorem \ref{thm:2} is essentially optimal under such a prescribed main term.

\begin{theorem}\label{thm:3}
Let $c_1$ and $c_2$ be two given positive numbers.
Suppose that $f(x)\ge 1$ is a monotone increasing function such that $\lim_{x\to\infty}f(x)=\infty$, $f(n^2)=O\big(f(n)\big)$ and $f(2n)=f(n)+O(1)$. Then there are two subsets $\mathcal{A}$ and $\mathcal{B}$ of natural numbers with 
$$
\mathcal{A}(x)=c_1x^{\alpha}+O\left(\frac{x^\alpha}{f(x)}\right) \quad \text{and} \quad 
\mathcal{B}(x)=c_2x^{1-\alpha}+O\left(\frac{x^{1-\alpha}}{f(x)}\right)
$$
such that for all sufficiently large $x$,
$$
\sum_{n\le x}d_{\mathcal{B},\mathcal{A}}(n)\ll x \log f(x).
$$
\end{theorem}

The two auxiliary conditions on $f$ in Theorem \ref{thm:3} are independent; see Example \ref{example:2}.

Theorem \ref{thm:2} and Theorem \ref{thm:3} give rather satisfactory answers to our initial problem.
We can remove the auxiliary conditions in Theorem \ref{thm:3} by an elementary argument to form the following clean result. The payoff is that we no longer control the error terms in the counting-function estimates for $\mathcal{A}(x)$ and $\mathcal{B}(x)$.

\begin{theorem}\label{thm:3-2}
Let $c_1$ and $c_2$ be two given positive numbers. Let \(f(x)>0\)  be a monotone increasing function with \(f(x)\to\infty\) as \(x\to\infty\). Then there exist $\mathcal{A}, \mathcal{B}\subset \N$ with $\mathcal{A}(x)\sim c_1x^{\alpha}$ and $\mathcal{B}(x)\sim c_2x^{1-\alpha}$ as $x\to\infty$ such that, for all sufficiently large $x$,
	\[
	\sum_{n\leq x}d_{\mathcal{B},\mathcal{A}}(n)\ll x f(x).
	\]
\end{theorem}

Our last theorem considers a single value of $d_{\mathcal{B},\mathcal{A}}(n)$.

	\begin{theorem}\label{thm:4}
Let $c>0$ be a given number and $\mathcal{A}$ a subset of $\N$ with $\mathcal{A}(x)\sim cx^\alpha$ as $x\to\infty$.
		Then  there exists a subset $\mathcal{B}$ with $\mathcal{B}(x)\asymp x^{1-\alpha}$ such that for $n\geq 2$,
		\[
        d_{\mathcal{B}, \mathcal{A}}(n)\ll \log n.
		\]
	\end{theorem}
The paper is organized as follows. In Section 2, discuss some examples and remarks which motivates our results. In Section 3, we prove our theorems. Finally, in Section 4 we provide a few open problems for further research.


\section{Examples and remarks}\label{sec:examples}

We record here the examples and auxiliary remarks that were used to motivate the main statements.

\begin{remark}\label{remark:upper_bound}
		 It is worth noting that the upper bound condition on $\mathcal{B}(x)$ in Proposition \ref{thm:lb:ABasy} is necessary, which  at first glance may seem a bit counterintuitive. We offer here an example to illustrate this.
         Let
					\begin{align*}
					\mathcal{A}=\big\{2n^2: n\in \mathbb{N}\big\} \quad \text{and} \quad
					\mathcal{B}=\cup_{k\geq 1}\big[2^{2^{k}}+1,2^{2^{k}+1}\big].
				\end{align*}
				Clearly, we have  $\mathcal{A}(x)= (x/2)^{1/2}+O(1)$. Let $x$ be a sufficiently large number.
                Then there exists some $\ell$ such that
                $$
                2^{2^\ell}<x\le 2^{2^{\ell+1}}.
                $$
                If $2^{2^\ell}<x<2^{2^\ell+1}$, then
                \begin{align*}
                \mathcal{B}(x)\geq x-2^{2^\ell}+2^{2^{\ell-1}}> 2^{2^{\ell-1}}=\frac{1}{\sqrt{2}}\sqrt{2^{2^{\ell}+1}}\gg\sqrt{x}.
                \end{align*}
                If $2^{2^{\ell}+1}\le x\le 2^{2^{\ell+1}}$, then
                \begin{align*}
                \mathcal{B}(x)\geq 2^{2^\ell}=\sqrt{2^{2^{\ell+1}}}\ge \sqrt{x}.
                \end{align*}
                Hence, we have $\mathcal{B}(x)\gg x^{1/2}$ for all large $x$.
                However, one notes easily that
                \begin{align*}
            \mathcal{B}\left(2^{2^{k}+1}\right)\ge 2^{2^{k}+1}-2^{2^{k}}=2^{2^{k}+1}/2,
                \end{align*}
                which shows that $\mathcal{B}(x)\gg x$ for infinitely many $x$. Hence, $\mathcal{B}(x)$ does not satisfy the upper bound of Proposition \ref{thm:lb:ABasy}.

                Now, let $a=2^{2^k+1}$ be a element of $\mathcal{A}$. 
                For any $b\in \mathcal{B}$ with $b>a$, it is easy to see
        $$
        b-a\ge 2^{2^{k+1}}+1- 2^{2^k+1}>2^{2^{k+1}}- 2^{2^k+1}.
        $$
        Therefore, for $x_k=2^{2^{k+1}}- 2^{2^k+1}$ we have
				\[
			\sum_{n\leq x_k}d_{\mathcal{B},\mathcal{A}}(n)\le \mathcal{A}\big(2^{2^k+1}\big)\mathcal{B}\big(2^{2^{k}+1}\big)\le 2^{2^{k-1}}\cdot 2^{2^{k}+1}\ll x_k^{1/4}x_k^{1/2}\ll x_k^{3/4}
				\]
                showing that $\sum_{n\leq x}d_{\mathcal{B},\mathcal{A}}(n)\ll x^{3/4}$ for infinitely many $x$.
\end{remark}
Now, a natural question could be asked. Whether the lower bound in Proposition \ref{thm:lb:ABasy} is tight. The following example shows that the bound given in Proposition \ref{thm:lb:ABasy} is optimal.

	\begin{example}\label{example:1}
	    Let $\mathcal{A}_0=\{\lfloor n/\alpha\rfloor: n\in \mathbb{N}\}$ and $\mathcal{B}_0=\mathbb{N}\setminus \mathcal{A}_0$. Then $\mathcal{A}_0(x)=\alpha x+O(1)$ and $\mathcal{B}_0(x)=(1-\alpha)x+O(1)$, where the implied constants in the big-O notation depend only on $\alpha$. Let
        $$
\mathcal{A}=\left\{\sum_{a}\varepsilon_a 2^{a}: a\in \mathcal{A}_0, \varepsilon_a\in\{0, 1\},\ \sum_{a}\varepsilon_a<\infty\right\}.
$$
and
$$
\mathcal{B}=\left\{\sum_{b}\varepsilon_b 2^{b}:b\in \mathcal{B}_0, \varepsilon_b\in\{0, 1\},\ \sum_{b}\varepsilon_b<\infty\right\}.
$$
By the constructions of $\mathcal{A}_0$ and $\mathcal{B}_0$, we have $\mathcal{A}(x)\asymp x^{\alpha}$ and $\mathcal{B}(x)\asymp x^{1-\alpha}$. By the unique binary representation of integers, we have $d_{\mathcal{B},\mathcal{A}}(n)\le 1$ for any $n\in \N$. Hence,
\[
\sum_{n\leq x}d_{\mathcal{B},\mathcal{A}}(n)\le x.
\] 
	\end{example}

The preceding example shows the sharpness of Proposition \ref{thm:lb:ABasy} without additional regularity assumptions on $\mathcal{A}(x)$.

A simple argument implies $f(x)\ll \log x$ from the condition $f(2n)=f(n)+O(1)$ (see the proof of Theorem \ref{thm:3}). It is also worth mentioning that there is no containment relationship between $f(2n)=f(n)+O(1)$ and $f(n^2)=O\big(f(n)\big)$. We provide examples below to illustrate this.

\begin{example}\label{example:2}
{\rm I.)} $f(n^2)=O\big(f(n)\big)$ does not imply $f(2n)=f(n)+O(1)$.

The function $f(x)=(\log x)^2$ is a simple example.

\medskip

{\rm II.)} $f(2n)=f(n)+O(1)$ does not imply $f(n^2)=O\big(f(n)\big)$.

\noindent
We define the function $f$ by recursions. Let $n_1$ be a large given integer.

\noindent
Let $f(x)=\sqrt{\log x}$ for $e\le x\le 2^{n_1}$. 

\noindent
For $1\le k\le n_1$, let $f(x)=f\big(2^{n_1+k-1}\big)+\log 2$ if $2^{n_1+k-1}<x\le 2^{n_1+k}$. 

\noindent
Let $f(x)=f\big(2^{2n_1}\big)$ for $2^{2n_1}\le x\le \exp\Big(f\big(2^{2n_1}\big)^2\Big):=n_2$.

\noindent
Let $f(x)=\sqrt{\log x}$ for $n_2\le x\le 2^{n_2}$. 

\noindent
For $1\le k\le n_2$, let $f(x)=f\big(2^{n_2+k-1}\big)+\log 2$ if $2^{n_2+k-1}<x\le 2^{n_2+k}$. 

\noindent
Let $f(x)=f\big(2^{2n_2}\big)$ for $2^{2n_2}\le x\le \exp\Big(f\big(2^{2n_2}\big)^2\Big):=n_3$.

\noindent
Suppose that $n_j$ is defined. We now define $n_{j+1}$.

\noindent
Let $f(x)=\sqrt{\log x}$ for $n_j\le x\le 2^{n_j}$. 

\noindent
For $1\le k\le n_j$, let $f(x)=f\big(2^{n_j+k-1}\big)+\log 2$ if $2^{n_j+k-1}<x\le 2^{n_j+k}$. 

\noindent
Let $f(x)=f\big(2^{2n_j}\big)$ for $2^{2n_j}\le x\le \exp\Big(f\big(2^{2n_j}\big)^2\Big):=n_{j+1}$.

By Continuing this process, we define the function $f$. Clearly, we have $f(n)\le f(2n)\le f(n)+\log 2$ by this definition. However, it is not hard to see
$$
\frac{f(2^{2n_j})}{f(2^{n_j})}\gg \sqrt{n_j}\to\infty,\ \text{as~} j\to\infty.
$$
\end{example}

\section{Probabilistic tools}\label{sec:probabilistic tools}
In this section we give a short summary about the probabilistic method, which plays an important role in some of our proofs. Interested reader can find more details in \cite{TV06}. In this paper we denote the probability of an event $E$ by $\mathbf{P}(E)$ and the expectation of a random variable $X$ by $\mathbb{E}(X)$. Let $\mathbf{1}_{\mathcal{A}}$ be the indicator random variable defined as
\[
\mathbf{1}_{\mathcal{A}}(n) = \begin{cases}
		1, & \text{if}\ n\in \mathcal{A}, \\
		0, & \text{if}\ n\notin \mathcal{A}.
	\end{cases}
\]
Then for a fixed set $\mathcal{A}$ and a randomly defined set $\mathcal{B}$ 
\[
d_{\mathcal{B},\mathcal{A}}(n) = \sum_{a\in \mathcal{A}}\mathbf{1}_{\mathcal{B}}(a+n)
\]
is a sum of independent indicator a random variables. 
Then we have
$$
\E\big(d_{\mathcal{B},\mathcal{A}}(n)\big)=\E\Big(\sum_{a\in \mathcal{A}}\mathbf{1}_{\mathcal{B}}(a+n)\Big)=\sum_{a\in\mathcal{A}}\bP\big((a+n)\in \mathcal{B}\big).
$$
Furthermore,
\[
\mathcal{A}(x) = 
\sum_{\substack{a \le x \\ a\in\mathbb{N}}}\mathbf{1}_{\mathcal{A}}(a).
\]
 Obviously,
$$
        \E\big(\mathcal{A}(x)\big)=\E\Big(\sum_{\substack{a\leq x \\ a\in \mathbb{N}}} \mathbf{1}_{\mathcal{A}}(a)\Big)=\sum_{a\le x}\bP\big(a\in \mathcal{A}\big).
        $$
We need the following inequality which plays a crucial role in our proofs; see \cite[p. 11]{TV06}.

\begin{proposition}[Chernoff's inequality] \label{Chernoff Bound}
					Let $X_1, X_2, \dots, X_n$ be independent boolean (i.e., 0-1 valued) random variables.
Define
\[
X = \sum_{i=1}^n X_i, \qquad \mu = \mathbb{E}(X).
\]
Then the following bounds hold:

\medskip

\noindent $(i)$ \textit{$($Upper tail$)$} For any $\delta > 0$,
\[
\bP\!\big(X \ge (1+\delta)\mu\big)
\le
\exp\!\left(-\frac{\delta^2}{2+\delta}\mu\right).
\]

\medskip

\noindent $(ii)$ \textit{$($Two-sided bound$)$} For any $0<\delta<1$,
\[
\bP\!\big(|X-\mu|\ge \delta\mu\big)
\le
2\exp\!\left(-\frac{\delta^2}{3}\mu\right).
\]
					
\end{proposition}

We will also use the Borel-Cantelli Lemma \cite[p. 3]{TV06}.

\begin{proposition}[Borel--Cantelli]\label{Borel}
					Let $(\Omega, \mathcal{F}, \mathbb{P})$ be a probability space, and let
					$\{E_n\}_{n=1}^\infty$ be a sequence of events in $\mathcal{F}$.

					If
					\[
					\sum_{n=1}^{\infty} \bP(E_n) < \infty,
					\]
					then with probability 1 at most a finite number of events $E_{n}$ can occur.
\end{proposition}

\section{Proofs}\label{sec:proofs}

\subsection{Proof of Proposition \ref{thm:lb:ABasy}}
\begin{proof}
        By the conditions on $\mathcal{A}$ and $\mathcal{B}$, there are positive constants $c_i~(1\le i\le 4)$ such that
$$
c_1x^{\alpha}<\mathcal{A}(x)<c_2x^{\alpha} \quad \text{and} \quad c_3x^{1-\alpha}<\mathcal{B}(x)<c_4x^{1-\alpha}.
$$
Let $c>0$ be a sufficiently small constant to be decided later. Clearly,
\begin{align*}
\sum_{n\leq x}d_{\mathcal{B},\mathcal{A}}(n)&\ge \sum_{\substack{a\le cx,\ cx<b\le x\\ a\in \mathcal{A},\ b\in\mathcal{B}}}1\\
&= \mathcal{A}(cx)\big(\mathcal{B}(x)-\mathcal{B}(cx)\big)\\
&\geq c_1(cx)^{\alpha}\big(c_3x^{1-\alpha}-c_4(cx)^{1-\alpha}\big)\\
&=c_1c^\alpha\big(c_3-c_4c^{1-\alpha}\big)x.
\end{align*}
The proposition follows from choosing $c=\frac{1}{2}\big(\frac{c_3}{c_4}\big)^{\frac{1}{1-\alpha}}$.
\end{proof}

\subsection{Proof of Theorem \ref{thm:2}}
\begin{proof}[Proof of Theorem \ref{thm:2}]
    For the first assertion, we have to show
    \[
	\sum_{n\leq x}d_{\mathcal{B},\mathcal{A}}(n)\gg x\log f(x)
	\]
    for any $f(x)\le x^\alpha$. Let $x$ be sufficiently large. Let $k$ be a positive integer. 
Suppose that 
$$
b\in \mathcal{B} \cap \Big(kx/2, (k+1)x/2\Big] \quad \text{and} \quad  a\in \mathcal{A}\cap \Big((k-1)x/2, kx/2\Big].
$$
Then we clearly have
$1\le b-a\le x$. Hence,
\begin{align}\label{eq-Thm:2-1}
\sum_{n\le x}d_{\mathcal{B},\mathcal{A}}(n)\ge \sum_{k=1}^{\infty}\left(\mathcal{B}\Big(\frac{(k+1)x}{2}\Big)-\mathcal{B}\Big(\frac{kx}{2}\Big)\right)\left(\mathcal{A}\Big(\frac{kx}{2}\Big)-\mathcal{A}\Big(\frac{(k-1)x}{2}\Big)\right).
\end{align}
By the assumed estimate for $\mathcal{A}(x)$, for $k\ge 3$ we have
\begin{align*}
\mathcal{A}\Big(\frac{kx}{2}\Big)-\mathcal{A}\Big(\frac{(k-1)x}{2}\Big)=\frac{c(kx)^\alpha}{2^\alpha}\left(1-\Big(1-\frac{1}{k}\Big)^\alpha\right)+O\left(\frac{(kx)^\alpha}{f(x)}\right).
\end{align*}
Note that $1-(1-y)^\alpha \ge \alpha y$ for any $0<y<1$. Then
\begin{align*}
\mathcal{A}\Big(\frac{kx}{2}\Big)-\mathcal{A}\Big(\frac{(k-1)x}{2}\Big)\ge \frac{c\alpha(kx)^\alpha}{2^\alpha k}
+O\left(\frac{(kx)^\alpha}{f(x)}\right).
\end{align*}
Let $C>0$ denote the implied constant in the error term above.
Let $K_x$ be the largest integer $k$ such that
\begin{align}\label{eq-Thm:2-2}
\frac{2^{\alpha+1} Ck}{c\alpha}\le f(x). 
\end{align}
Hence, for any $3\le k\le K_x$ we have
\begin{align}\label{eq-Thm:2-3}
\mathcal{A}\Big(\frac{kx}{2}\Big)-\mathcal{A}\Big(\frac{(k-1)x}{2}\Big)\ge \frac{1}{2}\frac{c\alpha(kx)^\alpha}{2^{\alpha} k}:= \frac{c_1 x^\alpha}{k^{1-\alpha}},\ \text{say},
\end{align}
where $c_1=\frac{c\alpha}{2^{1+\alpha}}$.
Combining (\ref{eq-Thm:2-1}) with (\ref{eq-Thm:2-3}) we get
\begin{align*}
\sum_{n\le x}d_{\mathcal{B},\mathcal{A}}(n)\ge c_1 x^\alpha\sum_{3\le k\le K_x}\left(\mathcal{B}\Big(\frac{(k+1)x}{2}\Big)-\mathcal{B}\Big(\frac{kx}{2}\Big)\right)\frac{1}{k^{1-\alpha}}.
\end{align*}

Using the Abel summation, we have
\begin{align*}
\sum_{n\le x}d_{\mathcal{B},\mathcal{A}}(n)\ge c_1 x^\alpha\sum_{4\le k\le K_x}\mathcal{B}\Big(\frac{kx}{2}\Big)\left(\frac{1}{(k-1)^{1-\alpha}}-\frac{1}{k^{1-\alpha}}\right)-\frac{c_1x^\alpha}{3^{1-\alpha}} \mathcal{B}\Big(\frac{3x}{2}\Big).
\end{align*}
Since $\mathcal{B}(x)\asymp x^{1-\alpha}$, from the above estimate we obtain 
\begin{align*}
\sum_{n\le x}d_{\mathcal{B},\mathcal{A}}(n) \gg x\sum_{4\le k\le K_x}k^{1-\alpha}\left(\frac{1}{(k-1)^{1-\alpha}}-\frac{1}{k^{1-\alpha}}\right)+O(x).
\end{align*}
Noting that for any $k\ge 4$ we have
\begin{align*}
k^{1-\alpha}\left(\frac{1}{(k-1)^{1-\alpha}}-\frac{1}{k^{1-\alpha}}\right)=\left(1+\frac{1}{k-1}\right)^\alpha-1\gg \frac{1}{k},
\end{align*}
where the implied constant depends only on $\alpha$. It then follows that
\begin{align}\label{first-eq-1}
\sum_{n\le x}d_{\mathcal{B},\mathcal{A}}(n) \gg  x\sum_{4\le k\le K_x}\frac{1}{k}\gg x\log K_x.
\end{align}
By (\ref{eq-Thm:2-2}) we know
\begin{align}\label{first-eq-2}
K_x= \left\lfloor \frac{c\alpha f(x)}{2^{\alpha+1}C}\right\rfloor\gg f(x),
\end{align}
from which the first assertion clearly follows.

We now turn to the proof of the second assertion.
From (\ref{first-eq-1}) and (\ref{first-eq-2}) we know
$$
\left|\Big\{ (b,a):a\in \mathcal{A}, b\in \mathcal{B}, 0<b-a<x, b\le xK_x \Big\} \right|\gg x\log f(x),
$$
provided that $x$ is large enough. Note that 
\begin{align*}
\left|\Big\{ (b,a)\in \mathcal{B}\times\mathcal{A}:0<b-a<x^{\alpha ^2}, b\le xK_x   \Big\} \right|\le \mathcal{B}(xK_x)x^{\alpha ^2}\ll (x^{1+\alpha })^{1-\alpha } x^{\alpha ^2}\ll x.
\end{align*}
Hence,
$$
\left|\Big\{ (b,a)\in \mathcal{B}\times\mathcal{A}:x^{\alpha ^2}\le b-a<x, b\le xK_x   \Big\}\right |\gg x\log f(x).
$$
It follows that there exists an $x^{\alpha ^2} \le n<x$ such that 
$$
d_{\mathcal{B},\mathcal{A}}(n)\gg \log f(x) \gg \log f(n),
$$
which completes the proof.
\end{proof}



\subsection{Proof of Theorem \ref{thm:3}}
\begin{proof}[Proof of Theorem \ref{thm:3}]
We first show that $f(x)\ll \log x$. By the condition 
$$
f(2n)=f(n)+O(1),
$$
there exists an absolute constant $c_f>0$ such that
\begin{align}\label{eq-constant}
    f(2n)\le f(n)+c_f.
\end{align}
Let $x\ge 2$ be a given number. Then there exists some integer $k$ such that $2^k\le x<2^{k+1}$. By (\ref{eq-constant}) we have
$$
f(x)\le f\big(2^{k+1}\big)\le f\big(2^k\big)+c_f\le \cdots\le f(2)+kc_f\le f(2)+c_f\frac{\log x}{\log 2}\ll \log x.
$$

For $0<\alpha<1$, we choose a fixed number $\lambda_\alpha$ satisfying that $1-\alpha<\lambda_\alpha<1$. Then $0<1-\lambda_\alpha<\alpha$. We then choose another fixed number $\gamma_\alpha$ such that 
\begin{align}\label{eq-lambda-1}
    \frac{1}{\alpha}<\gamma_\alpha<\frac{1}{1-\lambda_\alpha}.
\end{align}
Since $\gamma_\alpha \lambda_\alpha>\gamma_\alpha-1$ by (\ref{eq-lambda-1}), it follows  that there exists some $k_0$ such that
\begin{align}\label{eq-lambda-2}
(k+1)^{\gamma_\alpha}<k^{\gamma_\alpha}+2\gamma_\alpha k^{\gamma_\alpha-1} \quad \text{and} \quad k^{\gamma_\alpha \lambda_\alpha}\ge 2\gamma_\alpha k^{\gamma_\alpha-1},
\end{align}
provided that $k\ge k_0$.
We will define the sequence $n_{k_0}<n_{k_0+1}<\dots $ by recursions.  

First, let $n_{k_0}=\lceil k_0^{\gamma_\alpha} \rceil $. Suppose that we have defined $n_k~(k\ge k_0)$, and then the integer $n_{k+1}$ is defined to be the least integer $\ell$ such that
\begin{align}\label{eq-lambda-3}
\ell>n_{k}+n_k^{\lambda_\alpha}
\end{align}
and
\begin{align}\label{eq-ell-1}
f(n_k)n_k+2n_k\le f(\ell)\ell,
\end{align}
We now prove $n_k\ge k^{\gamma_\alpha}~(k\ge k_0)$ by induction on $k$. First of all, we have $n_{k_0}\ge k_{0}^{\gamma_\alpha}$ by the definition of $n_{k_0}$. Now suppose that $n_k\ge k^{\gamma_\alpha}$, then by (\ref{eq-lambda-3}) and (\ref{eq-lambda-2}) we have
$$
n_{k+1}>n_k+n_k^{\lambda_\alpha}\ge k^{\gamma_\alpha}+k^{\gamma_\alpha \lambda_\alpha}\ge k^{\gamma_\alpha}+2\gamma_\alpha k^{\gamma_\alpha-1}\ge (k+1)^{\gamma_\alpha}.
$$
Hence, we conclude that 
\begin{align}\label{eq-lambda-4}
k\le n_k^{\frac{1}{\gamma_\alpha}}.
\end{align}

Note that 
$$
f(n_k)n_k+2n_k\le f\left(n_k+2\Big\lceil \frac{n_k}{f(n_k)}\Big\rceil \right)\left(n_k+2\Big\lceil \frac{n_k}{f(n_k)}\Big\rceil\right),
$$
and 
\[
n_k+2\Big\lceil \frac{n_k}{f(n_k)}\Big\rceil >n_{k}+n_{k}^{\lambda_\alpha},
\]
provided that $k$ is sufficiently large. Thus, there is some $k_0^*\ge k_0$ such that
$$
n_k+n_k^{\lambda_\alpha}< n_{k+1}\le n_k+2\Big\lceil \frac{n_k}{f(n_k)}\Big\rceil
$$
for any $k\ge k_0^*$. Hence, for $k\ge k_0^*$ we have
\begin{align}\label{eq-5-21-3}
n_{k+1}=n_k+O\left(\frac{n_k}{f(n_k)}\right),
\end{align}
since $\lambda_\alpha<1$ and $f(n_k)\ll \log n_k$.
Let
$$
A=\bigcup_{k=k_0^*}^{\infty}\Big[f(n_k)n_k+1,\ f(n_k)n_k+c_1f(n_k)^{\alpha}n_k^{\alpha}-c_1f(n_{k-1})^{\alpha}n_{k-1}^{\alpha}\Big]
$$
and
$$
B=\bigcup_{k=k_0^*}^{\infty}\Big[f(n_k)n_k-c_2f(n_k)^{1-\alpha}n_k^{1-\alpha}+ c_2f(n_{k-1})^{1-\alpha}n_{k-1}^{1-\alpha},\ f(n_k)n_k\Big].
$$

Next, we verify the counting-function estimates for $\mathcal{A}(x)$ and $\mathcal{B}(x)$. Note first that
\begin{align}\label{eq-5-21-1}
\mathcal{A}\Big(f(n_k)n_k+c_1&f(n_k)^{\alpha}n_k^{\alpha}-c_1f(n_{k-1})^{\alpha}n_{k-1}^{\alpha}\Big)\nonumber\\
&=c_1\sum_{i=k_0}^{k} \Big(f(n_i)^{\alpha}n_i^{\alpha}-f(n_{i-1})^{\alpha}n_{i-1}^{\alpha}\Big)+O(k)\nonumber\\
&=c_1f(n_k)^{\alpha}n_k^{\alpha}+O\left(n_k^{1/\gamma_\alpha}\right),
\end{align}
where the last estimate uses (\ref{eq-lambda-4}). Similarly, we have
\begin{align}\label{eq-5-21-2}
\mathcal{A}\Big(f(n_{k+1})n_{k+1}+c_1&f(n_{k+1})^{\alpha}n_{k+1}^{\alpha}-c_1f(n_{k})^{\alpha}n_{k}^{\alpha}\Big)\nonumber\\
&=c_1\sum_{i=k_0}^{k+1} \Big(f(n_i)^{\alpha}n_i^{\alpha}-f(n_{i-1})^{\alpha}n_{i-1}^{\alpha}\Big)+O(k)\nonumber\\
&=c_1f(n_{k+1})^{\alpha}n_{k+1}^{\alpha}+O\left(n_k^{1/\gamma_\alpha}\right).
\end{align}
Using (\ref{eq-5-21-3}) and the relation $f(2n)=f(n)+O(1)$, we obtain
\begin{align}\label{eq-5-21-4}
&f(n_{k+1})^{\alpha}n_{k+1}^{\alpha}=\left( f(n_k)+O(1) \right) ^{\alpha } \left(  n_k+O\left( \frac{n_k}{f(n_k)} \right) \right) ^{\alpha}\nonumber \\ =&\left( f(n_k)^{\alpha }+O\left( f(n_k)^{\alpha -1} \right)\right) \left( n_k^{\alpha} +O\left(  n_k^{\alpha -1}\frac{n_k}{f(n_k)} \right) \right) =f(n_{k})^{\alpha}n_{k}^{\alpha}+O\left(\frac{n_k^\alpha}{f(n_k)^{1-\alpha}}\right).
\end{align}
Now, let $x$ be a sufficiently large number such that
\begin{align*}
    f(n_{k})n_{k}+c_1f(n_{k})^{\alpha}n_{k}^{\alpha}&-c_1f(n_{k-1})^{\alpha}n_{k-1}^{\alpha}\\
    &\le x\le f(n_{k+1})n_{k+1}+c_1f(n_{k+1})^{\alpha}n_{k+1}^{\alpha}-c_1f(n_{k})^{\alpha}n_{k}^{\alpha}.
\end{align*}
Then, by (\ref{eq-5-21-3}) and the relation $f(2n)=f(n)+O(1)$ again, we get
\begin{align}\label{eq-5-21-5}
x=f(n_{k})n_{k}+O(n_k).
\end{align}
Since $1/\gamma_\alpha<\alpha$, by (\ref{eq-5-21-1}), (\ref{eq-5-21-2}), (\ref{eq-5-21-4}), and (\ref{eq-5-21-5}) we conclude that
\begin{align*}
\mathcal{A}(x)=c_1n_k^\alpha f(n_k)^{\alpha}+O\left(\frac{n_k^\alpha}{f(n_k)^{1-\alpha}}\right)=c_1x^\alpha+O\left(\frac{x^\alpha}{f(n_k)}\right)=c_1x^\alpha+O\left(\frac{x^\alpha}{f(x)}\right),
\end{align*}
where the last step uses the condition $f(n^2)\ll f(n)$ of the theorem. Similar arguments give the corresponding estimate for $\mathcal{B}(x)$:
$$
\mathcal{B}(x)=c_2x^{1-\alpha}+O\left(\frac{x^{1-\alpha}}{f(x)}\right).
$$

Let $x$ be sufficiently large. It remains to show that
$$
\sum_{n\le x}d_{\mathcal{B}, \mathcal{A}}(n)\ll x\log f(x).
$$
Suppose that $n_{h-1}<x\le n_h$. If $a\in\mathcal{A}, b\in\mathcal{B}$ satisfy $a<b$ and $b>n_hf(n_h)$, then there exist $i>h$ and $j<i$ such that
$$
b\in \Big[f(n_i)n_i-c_2\big({f(n_i)n_i}\big)^{1-\alpha},\ f(n_i)n_i\Big]
$$
and
$$
a\in \Big[f(n_j)n_j,\ f(n_j)n_j+c_1\big(f(n_j)n_j\big)^{\alpha}\Big].
$$
By definition (\ref{eq-ell-1}), we have
\begin{align*}
f(n_i)n_i\ge f(n_{i-1})n_{i-1}+2n_{i-1}.
\end{align*}
Since $j\le i-1$ and $i-1\ge h$, we have
\begin{align*}
    b-a&\ge f(n_i)n_i-c_2\big({f(n_i)n_i}\big)^{1-\alpha}-f(n_j)n_j-c_1\big(f(n_j)n_j\big)^{\alpha}\\
    &\ge 2n_{i-1}-O\Big(f(n_i)^{1-\alpha}n_i^{1-\alpha}+f(n_i)^{\alpha}n_i^{\alpha}\Big)\\
    & >n_h.
\end{align*}
Therefore, it follows from the above analysis and the counting-function estimates for $\mathcal{A}(x), \mathcal{B}(x)$ that
\begin{align*}
   &\sum_{n\le x}d_{\mathcal{B}, \mathcal{A}}(n)\le \sum_{s\le f(n_{h-1})}\Big(\mathcal{A}(sn_h)-\mathcal{A}\big((s-1)n_h\big)\Big)\Big(\mathcal{B}\big((s+1)n_h\big)-\mathcal{B}\big((s-1)n_h\big)\Big)\\
   &\ll n_h\sum_{2\le s\le f(n_{h-1})}\Big(s^\alpha-(s-1)^\alpha +O\big( \frac{1}{f(n_{h-1})} \big) \Big)\Big((s+1)^{1-\alpha}-(s-1)^{1-\alpha} ++O\big( \frac{1}{f(n_{h-1})} \big)\Big).
\end{align*}
By the Lagrange mean value theorem, we have
\begin{align*}
\Big(s^\alpha-(s-1)^\alpha\Big)\Big((s+1)^{1-\alpha}-(s-1)^{1-\alpha}\Big)\le \alpha (s-1)^{\alpha-1}(1-\alpha) (s-1)^{-\alpha}
\ll \frac{1}{s-1}.
\end{align*}
Hence, 
\begin{align*}
&\sum_{n\le x}d_{\mathcal{B}, \mathcal{A}}(n)\ll n_h\sum_{2\le s\le f(n_{h-1})}\frac{1}{s-1}+\frac{n_h}{f(n_{h-1})}\sum_{2\le s\le f(n_{h-1})}(s-1)^{\alpha -1}+\frac{n_h}{f(n_{h-1})}\sum_{2\le s\le f(n_{h-1})}(s-1)^{-\alpha } 
\\ & \ll n_h\log f(n_h) +\frac{n_h}{f(n_{h-1})^{1-\alpha }} +\frac{n_h}{f(n_{h-1})^{\alpha }}\ll x \log f(x),
\end{align*}
where the last step comes from (\ref{eq-5-21-3}) and the condition $f(n^2)\ll f(n)$,
completing the proof of our theorem.
\end{proof}

\subsection{Proof of Theorem \ref{thm:3-2}}
The proof of Theorem \ref{thm:3-2} makes use of the following elementary lemma.

\begin{lemma}\label{lem-1}
    Let $f$ be any monotonically increasing  function with $f(x)>0$ and $f(x)\rightarrow\infty $ as $x\to\infty$. Then there exists an increasing monotonic function $g$ satisfying that 

    $1)$ $g(x)\le f(x)$ for any $x\ge 1$,

   $2)$ $g(n^2)\le 2g(n)$ for any positive integer $n$,

   $3)$ $g(2n)=g(n)+O(1)$ for any positive integer $n$,

   $4)$ $g(x)\rightarrow\infty$ as $x\to\infty$.
\end{lemma}

\begin{proof}
    Suppose that $f\left(2^{2k}\right)=y_k$ for any $k\in \mathbb{N}$. Then $y_k\rightarrow\infty$ as $k\to\infty$.

\noindent
Let $g(x)=0$ for $1\le x<4$.

\noindent
Let $g(x)=b_1:=y_1$ for $2^2\le x<2^{4}$. 

\noindent
Let $g(x)=b_2:=b_1+\delta_1$ for $2^4\le x<2^{8}$, where 
$$
\delta_1=\min\big\{y_1, y_2-b_1\big\}.
$$

\noindent
Let $g(x)=b_3:=b_2+\delta_2$ for $2^8\le x<2^{16}$, where 
$$
\delta_2=\min\big\{y_1, y_3-b_2\big\}.
$$

\noindent
Continue this process. Suppose that we have defined the values of $g$ up to $x<2^{2k}$. 

\noindent
Let
$g(x)=b_{k}:=b_{k-1}+\delta_{k-1}$ for $2^{2k}\le x<2^{2(k+1)}$,
where 
$$
\delta_{k-1}=\min\big\{y_1, y_{k}-b_{k-1}\big\}.
$$

\noindent
First, by the definition of $g$, we clearly have $g(x)\le f(x)$. Condition 1) is satisfied.

\noindent
Second, for any large $n$, suppose $2^{2k}\le n<2^{2(k+1)}$. Then $2^{2(k+1)}\le n^2<2^{2(k+2)}$, and
$$
g(n^2)=b_{k+1}\le b_k+y_1=b_k+b_1\le 2b_k.
$$
Condition 2) is satisfied.

\noindent
Third, for any large $n$, suppose $2^{2k}\le n<2^{2(k+1)}$. Then $2^{2k+1}\le 2n<2^{2(k+1)+1}$, and 
$$
g(2n)\le b_{k+1}\le b_k+y_1=g(n)+O(1).
$$
Condition 3) is satisfied.

To verify the fourth condition, we separate the discussion into two cases. If there are infinitely many $k$ such that $y_{k}-b_{k-1}\le y_1$, then for these $k$ we have
$$
b_{k}=b_{k-1}+\delta_{k-1}=b_{k-1}+y_{k}-b_{k-1}=y_k\rightarrow\infty,\ \text{as}~k\to\infty. 
$$
Since $g$ is an increasing monotonic function, we have $g(x)\rightarrow\infty$ as $x\to\infty$.

If there are only finitely many $k$ such that $y_{k}-b_{k-1}\le y_1$, then there exists some $k_0$ such that for all $k\ge k_0$ we have
$$
\delta_{k-1}=\min\big\{y_1, y_{k}-b_{k-1}\big\}=y_1.
$$
Therefore, for any positive integer $m$,
$$
b_{k_0+m}\ge b_{k_0}+my_1\rightarrow\infty,\ \text{as}~m\to\infty.
$$
Since $g$ is an increasing monotonic function, we have $g(x)\rightarrow\infty$ as $x\to\infty$. Condition 4) is satisfied.
\end{proof}

\begin{proof}[Proof of Theorem \ref{thm:3-2}]
 Let $g$ be the function defined in Lemma \ref{lem-1} above. Then $g$ satisfies the conditions in Theorem \ref{thm:3}. By Theorem \ref{thm:3} there are $\mathcal{A}$ and $\mathcal{B}$ with
 $$
 \mathcal{A}(x)=c_1x^\alpha+O\left(\frac{x^\alpha}{g(x)}\right)\sim c_1x^\alpha
 $$
 and
 $$
 \mathcal{B}(x)=c_2x^{1-\alpha}+O\left(\frac{x^{1-\alpha}}{g(x)}\right)\sim c_2x^{1-\alpha}
 $$
 such that
 $$
 \sum_{n\le x}d_{\mathcal{B}, \mathcal{A}}(n)\ll x \log g(x)\ll xg(x)\ll xf(x),
 $$
 completing the proof of our theorem.
\end{proof}

	\subsection{Proof of Theorem \ref{thm:4}}
	
	We will use the probabilistic method to prove Theorem \ref{thm:4}.

	Let $\mathcal{B}$ be the random set obtained by selecting each element of $\N$, independently, with
    \begin{align}\label{eq-random}
	\bP\big(b\in \mathcal{B}\big)=
	\begin{cases}
		0, & \text{if}\ \text{there~exist}~ a\in \mathcal{A}~ \text{and}~ 1\le k\leq \frac{1}{3c}a^{1-\alpha} \, \text{such that } b=a+k, \\
		\frac{1}{b^{\alpha}}, & \text{otherwise.}
	\end{cases}
	\end{align}
	Define \[\mathcal{B}_{\text{drop}}=\bigcup_{a\in \mathcal{A}}\left[a,a+\frac{1}{3c}a^{1-\alpha}\right].\]
	Then, we have the following results.

    \begin{lemma}\label{lemma:basic_properties}
	The following properties hold.

	$(1)$ For all sufficiently large $x$, the set 
    \(\mathcal{B}_{{\rm drop}}\) satisfies
        \[
		\mathcal{B}_{{\rm drop}}(x)\leq x/3.
		\]
		
	$(2)$ The expected counting function of \(\mathcal{B}\) satisfies
        \[
		\E\big(\mathcal{B}(x)\big) \asymp x^{1-\alpha}.
		\]
		
	$(3)$	For integer \(n\ge 2\),
		\[
		\E\big(d_{\mathcal{B},\mathcal{A}}(n)\big)\ll \log n.
		\]
\end{lemma}

\begin{proof}
	(1) By the definition of $\mathcal{B}_{{\rm drop}}$, we know that
    \[
    \mathcal{B}_{{\rm drop}}(x)\leq\sum_{\substack{a\in \mathcal{A} \\ \ a\leq x}}\frac{1}{3c}a^{1-\alpha}.
    \]
Partial summation yields
    \begin{align*}
         \sum_{\substack{a\in \mathcal{A} \\ a\leq x}}\frac{1}{3c}a^{1-\alpha}&=\frac{\mathcal{A}(x)}{3c}x^{1-\alpha} -(1-\alpha)\int_{1}^{x}\frac{\mathcal{A}(t)}{3ct^\alpha}{\rm d} t\\
         &=\left(\frac{1}{3}+o(1)-\frac{1-\alpha}{3}\right)x\\
         &\le \frac{x}{3}.
    \end{align*}
    Hence, 
    $
     \mathcal{B}_{{\rm drop}}(x)\leq \frac{x}{3}.
    $
	\medskip
	
	(2) By the definition of $\mathcal{B}$, clearly
        \[
		\E\big(\mathcal{B}(x)\big)\leq \sum_{b\leq x}\frac{1}{b^\alpha}\ll x^{1-\alpha}.
        \]
    Since $\mathcal{B}(x)\le x/3$ by (1), we obtain
    \[	
    \E\big(\mathcal{B}(x)\big)=\sum_{\substack{b\leq x \\ b\notin \mathcal{B}_{{\rm drop}}}}\frac{1}{b^\alpha}\geq \sum_{\frac{2}{3}x\leq b\leq x}\frac{1}{b^\alpha}\gg x^{1-\alpha}.
    \]
	
	(3) According to the definitions,
		\begin{align*}
			\E\big(d_{\mathcal{B},\mathcal{A}}(n)\big)&=\sum_{a\in \mathcal{A}} \bP\big(a+n\in  \mathcal{B}\big)
			=\sum_{a\in \mathcal{A}}\frac{1}{(a+n)^{\alpha}}\mathbf{1}_{ \overline{\mathcal{B}}_{{\rm drop}}}(a+n).
			\end{align*}
            For any $a\in \mathcal{A}$, if $a+n\notin \mathcal{B}_{{\rm drop}}$, then
            $$
            a+n>a+\frac{1}{3c}a^{1-\alpha},
            $$
            which implies
	$a\leq(3cn)^{\frac{1}{1-\alpha}}$. 
    It follows that
	\begin{align*}
    \E\big(d_{\mathcal{B},\mathcal{A}}(n)\big)=\sum_{\substack{a\in \mathcal{A} \\ a\leq(3cn)^{\frac{1}{1-\alpha}}}}\frac{1}{(a+n)^{\alpha}}\mathbf{1}_{\overline{\mathcal{B}}_{{\rm drop}}}(a+n)\le \sum_{\substack{a\in \mathcal{A}\\ a\leq(3cn)^{\frac{1}{1-\alpha}}}}\frac{1}{(a+n)^{\alpha}}.
	\end{align*}
			Since $\mathcal{A}(x)=\big(c+o(1)\big)x^{\alpha}$, using the partial summation again, we obtain
			\begin{align*}
				\E\big(d_{\mathcal{B},\mathcal{A}}(n)\big)&\le\sum_{\substack{a\in\mathcal{A} \\ a\le n}}\frac{1}{(a+n)^{\alpha}}+ \sum_{\substack{a\in \mathcal{A} \\ n<a\le (3cn)^{\frac{1}{1-\alpha}}}}\frac{1}{(a+n)^{\alpha}}\\
				&\leq\sum_{\substack{a\leq n \\ a\in \mathcal{A}}}\frac{1}{n^{\alpha}}+\alpha\int_{n}^{(3cn)^{\frac{1}{1-\alpha}}}\frac{\mathcal{A}(t)}{t^{1+\alpha}}{\rm d}t\\
                &\ll 1+ \alpha\int_{n}^{(3cn)^{\frac{1}{1-\alpha}}}\frac{1}{t}{\rm d}t\\
                &\ll \log n,
				\end{align*}
                where the implied constant depends on $c$ and $\alpha$.
\end{proof}
	
To continue the argument, we will employ the Chernoff inequality.
	




					



					

\begin{lemma}\label{lem:order_B(x)}
		For all sufficiently large $x$, the random sets $\mathcal{B}$ defined by (\ref{eq-random}) satisfies the estimate $\mathcal{B}(x)\asymp x^{1-\alpha} $ with probability 1.
	\end{lemma}
	\begin{proof}
		
		Let $\mu=\E\big(\mathcal{B}(x)\big)$. Then $\mu\asymp x^{1-\alpha}$ by part (2) of Lemma \ref{lemma:basic_properties}. By (ii) of Chernoff's inequality, we have
		\begin{align}\label{eq-chernoff-1}
		\bP\Big(\big|\mathcal{B}(x)-\mu\big|\geq \frac{1}{2}\mu\Big)\leq 2\exp\big(-c_0x^{1-\alpha}\big),
		\end{align}
        for some constant $c_0>0$.
		Thus, from (\ref{eq-chernoff-1}) we have
		\[
		\sum_{n =1}^{\infty}\bP\Big(|\mathcal{B}(n)-\mu|\geq \frac{1}{2}\mu\Big)\le 2\sum_{n=1}^{\infty}\exp\big(-c_0n^{1-\alpha}\big) <+\infty.
		\]
		The Borel--Cantelli Lemma implies that with probability 1 the inequality
        $$
        |\mathcal{B}(n)-\mu|\geq \frac{1}{2}\mu
        $$
        only happens for finitely many $n$. For any of such sets $\mathcal{B}$, we have
        \begin{align*}
        \frac{1}{2}\mu\le \mathcal{B}(n)\le \frac{3}{2}\mu, \quad \text{i.e.}, \quad \mathcal{B}(n)\asymp n^{1-\alpha}
        \end{align*}
        provided that $n$ is sufficiently large, leading to our desired result.
	\end{proof}

    We are ready to prove Theorem \ref{thm:4}.

	\begin{proof}[Proof of Theorem \ref{thm:4}]

		Let $v= \mathbb{E}\big(d_{\mathcal{B},\mathcal{A}}(n)\big)$. By part (3) of Lemma \ref{lemma:basic_properties}, we have
		$v \le c_1 \log n$ for some constant $c_1>0$. By (i) of Chernoff's inequality, for any $\delta > 0$,
		\[
		\bP\Big(d_{\mathcal{B},\mathcal{A}}(n) \ge (1+\delta)v\Big)
		\le \exp\!\left(-\delta v \frac{\delta}{2+\delta}\right).
		\]
		
		Choose $\delta$ such that
		\[
		\delta v = (2+4c_1)\log n.
		\]
		Since $v \le c_1 \log n$, we have
		\[
		\delta \ge 4 + \frac{2}{c_1},
		\]
		which implies
		\[
		\frac{\delta}{2+\delta}
		\ge \frac{4+\frac{2}{c_1}}{6+\frac{2}{c_1}}
		> \frac{2}{3}.
		\]
        Hence,
        $$
        \delta v\frac{\delta}{2+\delta}>2\log n\cdot \frac{2}{3}=\frac{4}{3}\log n.
        $$
		Therefore,
		\[
		\bP\Big(d_{\mathcal{B},\mathcal{A}}(n) \ge v + \delta v\Big)
		\le \exp\!\left(-\delta v \frac{\delta}{2+\delta}\right)
		\le \exp\!\left(-\frac{4}{3}\log n\right)
		= n^{-4/3},
		\]
		from which it follows that
		\[
		\sum_{n= 1}^{\infty}\bP\Big(d_{\mathcal{B},\mathcal{A}}(n) \ge (2+4c_1)\log n\Big)\le \sum_{n=1}^{\infty}n^{-\frac{4}{3}}<\infty.
		\]
		Hence, by the Borel--Cantelli lemma, with probability $1$  we have
		\begin{align}\label{eq-Borel-2}
		d_{\mathcal{B},\mathcal{A}}(n) < (2+4c_1)\log n\ll \log n
		\end{align}
		for all sufficiently large $n$.
		
		Combining (\ref{eq-Borel-2}) with Lemma \ref{lem:order_B(x)}, our theorem is true not only for a single set $\mathcal{B}$ but also for almost all such random sets.
	\end{proof}

\section{Conclusion}\label{sec:conclusion}

The results above show that the first moment of $d_{\mathcal{B},\mathcal{A}}(n)$ is sensitive to the regularity of the counting functions, while Theorem \ref{thm:4} gives a logarithmic point-wise upper bound for a suitable choice of $\mathcal{B}$. We end with two open problems.

Let $\log^{(1)}n=\log n$ and, recursively, $\log^{(k+1)}n=\log(\log^{(k)}n)$.
In view of Theorem \ref{thm:3} and Theorem \ref{thm:4}, we ask the following problem.

\begin{problem}\label{problem:1}
Let $k$ be a given positive integer. Does there exist a constant $c>0$ and two
subsets $\mathcal{A}, \mathcal{B}\subset \N$ with $\mathcal{A}(x)\sim cx^\alpha$ as $x\to\infty$ and $\mathcal{B}(x)\asymp x^{1-\alpha}$
such that
$$
d_{\mathcal{B}, \mathcal{A}}(n)\ll \log^{(k)} n?
$$
\end{problem}



Now, let us turn back to the problem of additive complements. Recall that $\mathscr{S}_k$ denotes the set of $k$-th powers.  von Doorn \cite{vD25} constructed an additive complement $\mathscr{W}$ of $\mathscr{S}_2$ that satisfies 
$$
\sum_{n\le x}r_{\mathscr{S}_2, \mathscr{W}}(n)\ll x;
$$
see also \cite[$\#$33]{Blo22}.
 Let $\mathcal{A}$ be a subset with $\mathcal{A}(x)=x^\alpha+O(1)$. Using von Doorn's construction, we can find an additive complement $\mathcal{B}$ of $\mathcal{A}$ such that
$$
\sum_{n\le x}r_{\mathcal{A}, \mathcal{B}}(n)\ll x.
$$



\section{Acknowledgments}

 P.~P.~P. was supported by the National Research, Development and Innovation Office NKFIH (Grant Nr. K146387 and Excellence program, Grant Nr. 153829). P.~P.~P.~was also partially supported by the Institute for Basic Science (IBS-R029-C4). S\'andor Kiss was supported by the NKFIH Grants No. K129335, K146387, KKP 144059. This author was also supported by the National Research, Development and Innovation Office NKFIH (Excellence program, Grant Nr. 153829. W. N. was supported by Shaanxi NSF (No. 2025JC-QYCX-002), Shaanxi Fundamental Science Research Project for Mathematics and Physics (No.25JSZ007) and the CSC programs. Csaba S\'andor was supported by the NKFIH Grants No. K129335, KKP 144059 and the Lend\"ulet ''Momentum'' program of the Hungarian Academy of Sciences (MTA).

\bibliographystyle{plain}

\end{document}